%% file: doccourse.tex
\documentclass[11pt,oneside]{amsart}
\usepackage{a4wide,mathrsfs,amsthm,amssymb,amsmath,url}
\input theoremstyle

\input commands
\title[Algebraic and model theoretic methods in constraint satisfaction]{Algebraic and model theoretic methods\\ in constraint satisfaction}
\author{Michael Pinsker}

\begin{document}
\begin{abstract}
This text is related to the tutorials I gave at the Banff International Research Station and within a ``Doc-course'' at Charles University Prague in the fall of 2014. It describes my current research and some of the most important open questions related to it.
\end{abstract}
\maketitle
\section{Overview}\label{subsub:brief}

A \emph{function clone} is a set $\C$ of finitary functions on a set $D$ which is closed under composition and which contains all projections. More formally,
\begin{itemize}
\item whenever $f\in \C$ is $n$-ary, and $g_1,\ldots,g_n\in\C$ are $m$-ary, then the $m$-ary function $f(g_1,\ldots,g_n)$ defined by
$$(x_1,\ldots,x_m)\mapsto f(g_1(x_1,\ldots,x_m),\ldots,g_n(x_1,\ldots,x_m))$$ is  an element of $\C$;
\item for all $1\leq k\leq n<\omega$, $\mathscr C$ contains the \emph{$k$-th $n$-ary projection $\pi^n_k\colon D^n\To D$}, uniquely defined by the equation $\pi^n_k(x_1,\ldots,x_n)=x_k$.
\end{itemize}
There are two main sources of function clones:
\begin{itemize}
\item the term operations of any algebra $\A$ form a clone, the \emph{term clone} of $\A$ (and in fact, every clone is of this form);
\item the set of all operations which \emph{preserve} a given relational structure $\Gamma$ form a clone, called the \emph{polymorphism clone} of $\Gamma$ (certain clones, the \emph{topologically closed clones}, are of this form).
\end{itemize}
The first source of function clones makes them an object of primary interest in universal algebra, since many properties of an algebra, such as its subalgebras and congruences, only depend on its term operations. The second source links them with relational structures, and in particular, as we will see, with certain questions in complexity theory. 

This topic of this text are function clones over a countably infinite set and their applications in complexity theory. While the investigation of \emph{all} such clones is not very promising since in the general setting, hardly any positive structural results could be expected (cf. for example~\cite{GoldsternPinsker}), research on clones which are ``sufficiently rich'' has proven extremely fruitful in recent years \cite{Topo-Birk, BP-reductsRamsey, OligoClone, MasulovicPech}. We are interested here in function clones which are rich in the sense that they contain a rather large permutation group: a permutation group on a  countably infinite set $D$ is called \emph{oligomorphic} iff its componentwise action on $D^n$ has only finitely many orbits, for all $n\geq 1$ (cf.~\cite{Oligo}). Function clones containing an oligomorphic permutation group, referred to as \emph{oligomorphic clones}~\cite{OligoClone}, have been shown to enjoy many properties of function clones on finite sets. For example, they satisfy a topological variant of Birkhoff's HSP theorem; moreover, they encode the complexity of certain computational problems, so-called \emph{constraint satisfaction problems (CSPs)}, and indeed have proven to be a valuable tool in the study of the complexity of such problems in what is called the \emph{algebraic approach} to CSPs. Oligomorphic function clones encode a much larger class of CSPs than function clones over finite sets~\cite{Bodirsky-HDR}, and yet many tools from the finite carry over to the oligomorphic setting.

Of particular importance to us will be oligomorphic clones which arise from \emph{homogeneous} relational structures in a finite language, and some of our methods rely on \emph{Ramsey theory} and on connections of Ramsey-type theorems with topological dynamics. While the original motivation for studying function clones comes from universal algebra, and later and independently from constraint satisfaction problems, their study therefore also involves tools and concepts from model theory, combinatorics, and topological dynamics.

\section{The State of the Art}

\noindent{\bf Birkhoff's theorem for oligomorphic clones.\:} I will start by recalling the finite version of Birkhoff's HSP theorem. An \emph{algebra} is a structure with a purely functional signature. The \emph{clone of an algebra $\A$} with signature $\tau$, denoted by $\Clo(\A)$, is the set of all functions with finite arity on the domain $A$ of $\A$ which can be written as $\tau$-terms over $\A$. More precisely, every abstract $\tau$-term $t$ naturally induces a finitary function $t^\A$ on $A$, and $\Clo(\A)$ consists precisely of the functions of this form.

Let $\A$, $\B$ be algebras of the same signature $\tau$. The assignment $\xi$ from $\Clo(\A)$ to $\Clo(\B)$ which sends every element $t^\A$ of $\Clo(\A)$ to $t^\B$ is a well-defined function if and only if for all $\tau$-terms $s,t$ we have that $s^\B=t^\B$ whenever $s^\A=t^\A$. In that case, $\xi$ is in fact a surjective \emph{clone homomorphism}, and we then call $\xi$ the \emph{natural homomorphism} from $\Clo(\A)$ onto $\Clo(\B)$. In general, a \emph{clone homomorphism} is a function $\sigma\colon \C\To \D$, where $\C,\D$ are clones (possibly acting on different base sets), which sends functions in $\C$ to functions of the same arity in $\D$, every projection in $\C$ to the corresponding projection in $\D$, and which preserves composition, i.e., $\sigma(f(g_1,\ldots,g_n))=\sigma(f)(\sigma(g_1),\ldots,\sigma(g_n))$ for all $n$-ary $f\in\C$ and all $m$-ary $g_1,\ldots,g_n\in\D$, for all $n,m\geq 1$ (cf.~\cite{BPP-reconstructing}).

When $\cC$ is a class of algebras with common signature $\tau$, then $\PPP(\cC)$ denotes the class of all products of algebras from $\cC$, $\PPPfin(\cC)$ denotes the class of all \emph{finite} products of algebras from $\cC$,
$\SSS(\cC)$ denotes the class of all subalgebras of algebras from $\cC$, and $\HHH(\cC)$ denotes the class of all 
homomorphic images of algebras from $\cC$ (when defining these operators, we consider algebras up to isomorphism). A \emph{pseudovariety} is a class $\cV$ of algebras
of the same signature 
such that $\cV = \HHH(\cV) = \SSS(\cV) = \PPPfin(\cV)$, i.e., a class closed under homomorphic images, subalgebras, and finite products; the pseudovariety \emph{generated} by a class of algebras $\cC$ (or by a single algebra $\A$) is the smallest pseudovariety that contains $\cC$ (contains $\A$, respectively). 
For \emph{finite} algebras, Birkhoff's HSP theorem takes the following form (see Exercise 11.5 in combination with the proof of Lemma 11.8 in~\cite{BS}).

\begin{thm}[Birkhoff~\cite{Bir-On-the-structure}]\label{thm:birk}
Let $\A,\B$ be finite algebras with the same signature. 
Then the following three statements are equivalent. 
\begin{enumerate}
\item The natural homomorphism from $\Clo(\A)$ onto $\Clo(\B)$ exists. 
\item $\B \in \HSPfin(\A)$.
\item $\B$ is contained in the pseudovariety generated by $\A$.
\end{enumerate}
\end{thm}
When $\A$ and $\B$ are of arbitrary cardinality, then the equivalence of $(2)$ and $(3)$ still holds; however, if one wants to maintain equivalence with item~(1), then another version of Birkhoff's theorem states that one has to replace
finite powers by arbitrary powers in the second item,
that is, one has to replace $\HSPfin(\A)$
by $\HSP(\A)$; the third item has to be adapted using the notion of a \emph{variety} of algebras, i.e., a class of algebras of common signature 
closed under the operators $\HHH$, $\SSS$ and $\PPP$.

It recently turned out that one can prove a similar theorem with finite powers for algebras $\A$ on a countably infinite domain whose clone $\Clo(\A)$ is oligomorphic -- we call such algebras \emph{oligomorphic} as well. To this end, one has to see function clones not only as algebraic, but also as topological objects.  On any set $D$, there is a largest function clone $\mathscr O_D$: the clone of all finitary operations on $D$. The ``function space'' $\mathscr O_D$ carries a natural topology, namely the topology of pointwise convergence, with respect to which the composition of functions is continuous. A basis of open sets of this topology 
is given by the sets of the form
$$\{ f \colon D^k\To D \; | \; 
f(a^1_1,\dots,a^1_k) = a^1_0, \dots,f(a^n_1,\dots,a^n_k) = a_0^n \} \;.$$
In fact, similarly to the Baire space ${\mathbb N}^{\mathbb N}$, 
$\mathscr O_D$ then becomes a Polish space (cf.~for example~\cite{BPP-reconstructing}). As a subset of $\mathscr O_D$, every function clone $\C$ on $D$ inherits this topology, and hence carries a topological structure in addition to its algebraic structure given by the equations which hold in $\C$. We denote the topological closure of a function clone $\C$ in $\mathscr O_D$ by $\overline{\C}$.

It is not hard to see that all algebras in the pseudovariety generated by an oligomorphic algebra are again oligomorphic. The following is the topological variant of Birkhoff's theorem for oligomorphic algebras.
 
 \begin{thm}[Bodirsky and Pinsker~\cite{Topo-Birk}]\label{thm:topobirk}
 Let $\A,\B$ be oligomorphic or finite algebras with the same signature. 
Then the following three statements are equivalent. 
\begin{enumerate}
\item The natural homomorphism from $\overline{\Clo(\A)}$ onto $\overline{\Clo(\B)}$ exists and is continuous. 
\item $\B \in \HSPfin(\A)$.
\item $\B$ is contained in the pseudovariety generated by $\A$.
\end{enumerate}
\end{thm}

Note that Theorem~\ref{thm:birk} really is a special case of Theorem~\ref{thm:topobirk}, since the topology of any function clone on a finite set is discrete, and hence the natural homomorphism from the clone of a finite algebra to that of another algebra is always continuous.\bigskip

\noindent{\bf Applications to constraint satisfaction problems.\:} Let us now turn to applications of oligomorphic function clones to computational complexity problems. Every relational structure $\Gamma$ in a finite language defines a computational problem, called the \emph{constraint satisfaction problem} of $\Gamma$ and denoted by $\Csp(\Gamma)$, as follows: input of the problem is a primitive positive sentence $\phi$ in the language for $\Gamma$, i.e., a sentence of the form $\exists x_1,\dots,x_n (\phi_1 \wedge \dots \wedge \phi_m)$ where $\phi_1,\dots,\phi_m$ are atomic formulas; the problem is to decide whether or not $\phi$ holds in $\Gamma$. An instance of this problem therefore asks about the existence of elements of $\Gamma$ satisfying a given conjunction of atomic conditions. The structure $\Gamma$ is called the \emph{template} of the problem, and can be finite or infinite. We will later see how infinite templates can model natural computational problems, and refer also to~\cite{Bodirsky-HDR} for an abundance of examples. We remark that $\Csp(\Gamma)$ is often presented in the form of a \emph{homomorphism problem}, which is easily seen to be equivalent: in this formulation, the input is a \emph{finite} structure $\Omega$ in the language of $\Gamma$ (which can still be finite or infinite), and the question is whether or not there exists a homomorphism from  $\Omega$ into $\Gamma$.

To every relational structure $\Gamma$, one can assign a function clone on the domain of $\Gamma$ as follows. A \emph{polymorphism} of a structure $\Gamma$ is a homomorphism from $\Gamma^k$ to $\Gamma$ for some finite $k\geq 1$; the \emph{polymorphism clone} $\Pol(\Gamma)$ of $\Gamma$ is the set of all polymorphisms of $\Gamma$. It is easy to see that $\Pol(\Gamma)$ is a function clone which is closed in the pointwise convergence topology described above, and in fact, the closed function clones are precisely the function clones of the form $\Pol(\Gamma)$ for a relational structure $\Gamma$.

For finite relational structures $\Gamma$, the complexity of $\Csp(\Gamma)$ depends, up to polynomial time, only on $\Pol(\Gamma)$ (cf.~\cite{JBK,JeavonsAlgebra}); this fact is the basis of the approach to constraint satisfaction via clones. The same is true for structures with oligomorphic polymorphism clones~\cite{BodirskyNesetrilJLC}. But which structures have oligomorphic polymorphism clones? The answer can be found in a classical theorem of model theory, the theorem of Engeler, Svenonius, and Ryll-Nardzewski (see e.g.~the textbook~\cite{HodgesLong}). A countable structure $\Gamma$ is called 
\emph{$\omega$-categorical} iff all countable models of the first-order theory of $\Gamma$ are isomorphic to $\Gamma$. Now the theorem states that  the automorphism group $\Aut(\Gamma)$ of a countable structure $\Gamma$ is oligomorphic if and only if $\Gamma$ is $\omega$-categorical. It follows that the polymorphism clone of a countable structure $\Gamma$ is oligomorphic if and only if $\Gamma$ is $\omega$-categorical.

\begin{thm}[Bodirsky and Ne\v{s}et\v{r}il~\cite{BodirskyNesetrilJLC}]\label{thm:csp:functional}
Let $\Gamma, \Gamma'$ be $\omega$-categorical structures in finite relational languages which have the same domain. If $\Pol(\Gamma)\subseteq \Pol(\Gamma')$, then $\Csp(\Gamma')$ is polynomial-time reducible to $\Csp(\Gamma)$.
\end{thm}

As a consequence, for $\omega$-categorical structures $\Gamma$ the complexity of their CSP is still up to polynomial time encoded in their polymorphism clone, i.e., if $\pol(\Gamma')=\pol(\Gamma)$, then $\Csp(\Gamma)$ and $\Csp(\Gamma')$ are polynomial-time equivalent. 

The theory of the algebraic approach to CSPs goes much further, which brings us back to Birkhoff's HSP theorem. For a structure $\Gamma$, we call any algebra on the domain on $\Gamma$ whose functions are precisely the elements of $\Pol(\Gamma)$ indexed in some arbitrary way a  \emph{polymorphism algebra} of $\Gamma$. It can be shown that if $\Gamma$ and $\Gamma'$ are finite structures in a finite relational language, and a polymorphism algebra $\B$ of $\Gamma'$ is contained in the pseudovariety of a polymorphism algebra $\A$ of $\Gamma$, then $\Csp(\Gamma')$ is polynomial-time reducible to $\Csp(\Gamma)$~\cite{JBK,JeavonsAlgebra}. By Birkhoff's theorem, this is the case iff the natural homomorphism from $\Clo(\A)$ onto $\Clo(\B)$ exists. One then sees that this is the case if and only if there exists a surjective clone homomorphism from $\Pol(\Gamma)$ onto $\Pol(\Gamma')$. Hence, the complexity of the CSP of a finite relational structure $\Gamma$ only depends on the abstract structure of the clone $\pol(\Gamma)$. Similarly to abstract groups, \emph{abstract clones} can be formalized as multi-sorted algebras equipped with composition operations as well as with constant symbols for the projections (cf.~, for example,~\cite{GoldsternPinsker} or~\cite{BPP-reconstructing}), but one can avoid this technicality: in practice, it is enough to know that abstract clones simply encode the equations which hold between its functions, or more precisely, it is enough to know that clone homomorphisms as defined above are precisely the structure preserving maps between those objects.

Using the topological generalization of Birkhoff's theorem, one can show the following for $\omega$-categorical structures (this is a simplified version; for a stronger formulation see~\cite{Topo-Birk}).

\begin{thm}[Bodirsky and Pinsker~\cite{Topo-Birk}]\label{thm:csp:topological}
Let $\Gamma, \Gamma'$ be finite or $\omega$-categorical structures in a finite relational language. If there exists a surjective continuous clone homomorphism from $\Pol(\Gamma)$ onto $\Pol(\Gamma')$, then $\Csp(\Gamma')$ is polynomial-time reducible to $\Csp(\Gamma)$.
\end{thm}

In particular, for $\omega$-categorical structures the complexity of their CSP is still up to polynomial time encoded in their polymorphism clone, seen as an abstract clone together with the topology on the functions. In analogy to topological groups, we call such objects \emph{topological clones}~\cite{BPP-reconstructing}. More precisely, if two $\omega$-categorical structures $\Gamma, \Gamma'$ have polymorphism clones which are isomorphic as topological clones (i.e., via a bijection which is a clone homomorphism, whose inverse is a clone homomorphism, and which is a homeomorphism), then their CSPs are polynomial-time equivalent.

A class of $\omega$-categorical structures for which the CSP is of particular interest are structures with a first-order definition in a homogeneous structure in a finite language. A structure $\Delta$ is called \emph{homogeneous} iff any isomorphism between finitely generated substructures of $\Delta$ extends to an automorphism of $\Delta$ (some authors call
this notion \emph{ultrahomogeneity} to distinguish it from related concepts
of homogeneity). A \emph{reduct} of a structure $\Delta$ is a relational structure on the same domain each of whose relations can be defined in $\Delta$ by a first-order formula without parameters. Countable homogeneous structures in a finite relational language are $\omega$-categorical, and reducts of $\omega$-categorical structures are $\omega$-categorical as well, and hence fall into our context (cf.~the textbook~\cite{HodgesLong}).

When $\Gamma$ is the reduct of a homogeneous structure in a finite language, then $\Csp(\Gamma)$ models a certain type of problem about finitely generated structures, as we will outline in the following. Homogeneous structures can be seen as generic objects, called \emph{Fra\"{i}ss\'{e} limits}, representing so-called \emph{Fra\"{i}ss\'{e} classes} of finitely generated structures. A \emph{Fra\"{i}ss\'{e} class} is a class ${\mathcal C}$ of finitely generated  structures in a fixed countable language closed under isomorphism and induced substructures which satisfies the {joint embedding property}, i.e., for all $\Omega_0, \Omega_1\in{\mathcal C}$ there is $\Omega_2\in {\mathcal C}$ such that $\Omega_0, \Omega_1$ embed into $\Omega_2$,  and the \emph{amalgamation property}, i.e., for any three structures $\Omega_0, \Omega_1, \Omega_2$ in ${\mathcal C}$  and embeddings $e:\Omega_0\To\Omega_1$ and $f:\Omega_0\To\Omega_2$ there exists $\Omega_3\in{\mathcal C}$ and embeddings $e':\Omega_1\To\Omega_3$ and $f':\Omega_2\To\Omega_3$ such that $e'\circ e=f'\circ f$. For any Fra\"{i}ss\'{e} class ${\mathcal C}$ there exists an up to isomorphism unique homogeneous structure $\Delta_{\mathcal C}$, called the Fra\"{i}ss\'{e} limit of ${\mathcal C}$, whose \emph{age}, i.e., the class of its finitely generated  substructures up to isomorphism, equals ${\mathcal C}$. Conversely, the age of any homogeneous structure in a countable language is a Fra\"{i}ss\'{e} class.

For example, the \emph{random graph} $G=(V,E)$ is the Fra\"{i}ss\'{e} limit of the class of finite undirected graphs without loops, and similarly there exist a \emph{random partial order}, a \emph{random tournament}, \emph{random hypergraphs}, a \emph{random digraph}, and so forth. Let us stick to the first example for a moment and let us define a class of computational problems about finite graphs as follows. Call quantifier-free formulas in the language of graphs \emph{graph formulas}. A graph formula $\Phi(x_1,\ldots,x_m)$ is \emph{satisfiable in a graph} iff there exists a graph $H$ and an $m$-tuple $a$ of elements in $H$ such that $\Phi(a)$ holds in $H$. Now let $\Psi = \{\psi_1,\dots,\psi_n\}$ be a finite set of graph formulas.
Then $\Psi$ gives
rise to the following computational problem.

\cproblem{Graph-SAT$(\Psi)$}
{A set of variables $W$ and a graph formula of the form
$\Phi = \phi_1 \wedge \dots \wedge \phi_l$ where each $\phi_i$ for $1 \leq i \leq l$ is obtained from one of the formulas $\psi$ in $\Psi$ by substituting the variables from $\psi$ by variables from $W$.}
{Is $\Phi$ satisfiable in a graph?}

In words, an instance of Graph-SAT$(\Psi)$ asks whether there exists a (finite) graph with satisfies a conjunction of properties; which properties can appear is restricted by the fixed set of graph formulas $\Psi$. Therefore, the computational complexity increases with $\Psi$ in the sense that $\Psi\subseteq\Psi'$, then any algorithm for Graph-SAT$(\Psi')$ solves Graph-SAT$(\Psi)$. It is easy to see that each problem Graph-SAT$(\Psi)$ is in NP, i.e., solvable in nondeterministic polynomial time. 
 
The connection with CSPs is that every problem Graph-SAT($\Psi$) can be translated into $\Csp(\Gamma_\Psi)$ for a finite language reduct $\Gamma_\Psi$ of the random graph $G=(V,E)$ and vice-versa. For one direction, let $\Psi=\{\psi_1,\ldots,\psi_n\}$ be a set of graph formulas. To this set we assign a reduct $\Gamma_\Psi$ of $G$ which has for each $\psi_i$ a relation $R_i$ consisting of those tuples of elements of $V$ that satisfy $\psi_i$ (where the arity of $R_i$ is given by the number of distinct variables that occur in $\psi_i$). One readily sees that any algorithm for Graph-SAT($\Psi$) can be adapted to $\Csp(\Gamma_\Psi)$ and vice-versa, and so the problems are essentially the same. For the other direction, if $\Gamma$ is a reduct of $G$ in a finite language, then each of its relations is defined by a first-order formula over $G$, and indeed even by a quantifier-free first-order formula (i.e., a graph formula), since homogeneity and $\omega$-categoricity  imply quantifier elimination. Let $\Psi_\Gamma$ the set of those graph formulas. Again, one easily checks that Graph-SAT($\Psi_\Gamma$)  and $\Csp(\Gamma)$ are basically the same problem~\cite{BodPin-Schaefer-both}.

Now let ${\mathcal C}$ be an arbitrary Fra\"{i}ss\'{e} class of finitely generated structures in a finite language. As in the case of graphs, for a finite set $\Psi$ of quantifier-free first-order formulas in the language of ${\mathcal C}$, we can define the following computational problem.

\cproblem{${\mathcal C}$-SAT$(\Psi)$}
{A set of variables $W$ and a formula of the form
$\Phi = \phi_1 \wedge \dots \wedge \phi_l$ where each $\phi_i$ for $1 \leq i \leq l$ is obtained from one of the formulas $\psi$ in $\Psi$ by substituting the variables from $\psi$ by variables from $W$.}
{Is $\Phi$ satisfiable in a structure in ${\mathcal C}$?}

As before, each problem ${\mathcal C}$-SAT$(\Psi)$ is equivalent to $\Csp(\Gamma_\Psi)$ for an appropriate reduct $\Gamma_\Psi$ of the Fra\"{i}ss\'{e} limit $\Delta_{\mathcal C}$ of ${\mathcal C}$ and vice-versa. Hence, classifying the complexity of the problems ${\mathcal C}$-SAT$(\Psi)$ and classifying the complexity of the constraint satisfaction problems of reducts of $\Delta_{\mathcal C}$ is one and the same thing. Complete classifications have been obtained so far for the following countable homogeneous structures.
\begin{itemize}
\item the empty structure $(\mathbb N,=)$~\cite{ecsps};
\item the order of the rationals $(\mathbb Q, \leq)$~\cite{tcsps-journal};
\item the random graph $G=(V,E)$~\cite{BodPin-Schaefer-both}. 
\end{itemize}

In each of the three cases the classifications resulted in dichotomies: the CSPs of the reducts turned out to be either NP-complete or in P (i.e., solvable in polynomial time, which we will henceforth refer to as \emph{tractable}). While there exist CSPs of homogeneous digraphs which are undecidable~\cite{BodPinTsa}, the following representability condition for a Fra\"{i}ss\'{e} class $\mathcal C$, arguably reasonable for the most interesting computational problems, forces ${\mathcal C}$-SAT problems to be in NP, and could possibly imply a general dichotomy. Let $\tau$ be a finite relational signature. A class ${\mathcal C}$ of finite $\tau$-structures is called \emph{finitely bounded} iff there exists a finite set of finite $\tau$-structures $\F$ such that ${\mathcal C}$ consists precisely of those finite $\tau$-structures which do not embed any element of $\F$. A relational structure is called finitely bounded iff its age is finitely bounded. When $\Gamma$ is a finite language reduct of a finitely bounded homogeneous structure $\Delta$, then $\Csp(\Gamma)$ is easily seen to be in NP. We conjecture the following.

\begin{conj}\label{conj:main}
Let $\Delta$ be  a finitely bounded homogeneous structure, and let $\Gamma$ be a finite language reduct of $\Delta$. Then $\Csp(\Gamma)$ is either in P or NP-complete.
\end{conj}

The advantage of translating ${\mathcal C}$-SAT problems into constraint satisfaction problems of reducts of the Fra\"iss\'{e} limit $\Delta_{\mathcal C}$ is that it allows for the algebraic approach via clones, as we will now outline.

Firstly, by Theorem~\ref{thm:csp:functional} we have that if $\Gamma, \Gamma'$ are $\omega$-categorical structures on the same domain and $\pol(\Gamma)\supseteq \pol(\Gamma')$, then $\Csp(\Gamma)$ has a polynomial-time reduction to $\Csp(\Gamma')$. On a theoretical level, this implies that when one wants to classify the complexity of the CSPs of all reducts of $\Delta_{\mathcal C}$, it suffices to consider all polymorphism clones of reducts of $\Delta_{\mathcal C}$ -- reducts with equal polymorphism clones are polynomial-time equivalent. Those clones are precisely the closed function clones which contain the automorphism group $\Aut(\Delta_{\mathcal C})$ of $\Delta_{\mathcal C}$. Moreover, if a closed function clone corresponds to a tractable (i.e., polynomial-time solvable) CSP, then so do all closed function clones containing it; if it corresponds to a NP-hard CSP, then so do all closed function clones above $\Aut(\Delta_{\mathcal C})$ contained in it. In the case of an existing dichotomy, one thus has to find the border between tractability and NP-hardness in the lattice of closed function clones containing $\Aut(\Delta_{\mathcal C})$. On a practical level, this implies that if the CSP of a reduct $\Gamma$ is in P, then this is witnessed by the presence of certain functions in $\Pol(\Gamma)$. And indeed, the presence of polymorphisms with certain properties have been successfully translated into algorithms in the classifications above -- see~\cite{BP-reductsRamsey}.

The second use of polymorphism clones is that Theorem~\ref{thm:csp:topological} allows us to compare CSPs on different domains, resulting in a tool both for showing tractability as well as for showing hardness. As for the latter, it turns out to be convenient to show NP-hardness of $\Csp(\Gamma)$ by exposing a continuous clone homomorphism from $\Pol(\Gamma)$ onto $\Pol(\Gamma')$, where $\Gamma'$ is an $\omega$-categorical or finite structure with a hard CSP~\cite{Topo-Birk}. In practice, $\Gamma'$ will generally be finite; in fact, $\Gamma'$ will often be any structure on a two-element set with a \emph{trivial} polymorphism clone, i.e., the clone ${\bf 1}$ of all projections on a two-element set (which is the polymorphism clone of NP-complete structures). Since ${\bf 1}$ is isomorphic to the smallest function clone (i.e., the clone of projections) on any finite set with at least two elements, and on finite domain smaller polymorphism clones correspond to harder CSPs, a continuous clone homomorphism to ${\bf 1}$ is in a sense the strongest finite reason for NP-hardness (more precisely, it implies that $\Gamma$ \emph{pp-interprets} all finite structures~\cite{Topo-Birk}; see also~\cite{BPP-projective}).

\begin{nota}
 We write ${\bf 1}$ for the clone of all projections on a two-element set.
\end{nota}

It is an open conjecture, and indeed the main conjecture for CSPs of finite structures known as the \emph{tractability conjecture}, that under a cosmetic assumption (the assumption of having an \emph{idempotent} polymorphism clone, see Section~\ref{sect:innovative}) on a finite structure $\Gamma$, the CSP for $\Gamma$ is NP-complete if there exists a clone homomorphism from $\pol(\Gamma)$ onto $\bf 1$, and in P otherwise. Clearly, if there exists such a homomorphism, then $\Csp(\Gamma)$ is NP-complete; the open part is the other direction. Note that if there is no such homomorphism, then this is witnessed by equations which hold in $\pol(\Gamma)$ but which cannot be satisfied in $\bf 1$. Numerous equations have been translated into algorithms, and indeed every non-trivial set of equations of an idempotent clone translates into an algorithm if one believes in the tractability conjecture. In the $\omega$-categorical setting, we cannot purely rely on equations, but need to take into account the topology on the functions -- at least with what we know today. Recently, research has been conducted investigating the role of this topology~\cite{BPP-projective}, and 
about how to show tractability in case there exists no continuous homomorphism to $\bf 1$.\bigskip

\noindent{\bf A Ramsey-theoretic method.\:} The behavior of polymorphism clones of reducts of homogeneous structures $\Delta$ in a finite relational language seems to be particularly close to that of function clones on finite sets when $\Delta$ satisfies a particular combinatorial property. Let $\tau$ be a relational signature. We say that a class ${\mathcal C}$ of finite $\tau$-structures is a \emph{Ramsey class} (in the sense of~\cite{NesetrilSurvey}) iff for all $\Omega_0,\Omega_1\in{\mathcal C}$ there exists $\Omega_2\in{\mathcal C}$ such that for all colorings of the copies of $\Omega_0$ in $\Omega_2$ with two colors there exists an isomorphic copy $\Omega_1'$ of $\Omega_1$ in $\Omega_2$ such that all copies of $\Omega_0$ in $\Omega_1'$ have the same color.
  A relational structure is called \emph{Ramsey} iff its age is a Ramsey class. When ${\mathcal C}$ is a relational Fra\"{i}ss\'{e} class and $\Delta_{\mathcal C}$ its Fra\"{i}ss\'{e} limit, then it is equivalent to call $\Delta_{\mathcal C}$ Ramsey iff for all $\Omega_0,\Omega_1\in{\mathcal C}$ and all colorings of the copies of $\Omega_0$ in $\Delta_{\mathcal C}$ there exists a copy of $\Omega_1$ in $\Delta_{\mathcal C}$ on which the coloring is constant. Examples of  Fra\"{i}ss\'{e} classes which are Ramsey classes are~the class of finite ordered undirected graphs, 
 the class of finite linear orders, and the class of finite partial orders with a linear extension~\cite{NesetrilRoedlPartite,NesetrilRoedlOrderedStructures,AbramsonHarrington}.

We remark that, for example, neither the random graph nor the random partial order are Ramsey, and thus seem to fall out of this framework. However, they are themselves reducts of homogeneous Ramsey structures, namely the random ordered graph (i.e., the Fra\"{i}ss\'{e} limit of the class of finite ordered undirected graphs) and the random partial order with a random linear extension (i.e., the Fra\"{i}ss\'{e} limit of the class of finite partial orders with a linear extension).

The Ramsey property can be exploited as follows. Let $\Xi$ be a structure. The \emph{type} of a tuple $b=(b^1,\ldots,b^n)$ of elements of $\Xi$, denoted by $\typ(b)$, is the set of first-order  formulas $\phi(x_1,\ldots,x_n)$ such that $\phi(b^1,\ldots,b^n)$ holds in $\Xi$. Now let $\Xi_1,\ldots,\Xi_m$ be structures. For an element $a$ of the product $\Xi_1\mult\cdots\mult\Xi_m$ and $1\leq i\leq m$, we write
    $a_i$ for the $i$-th coordinate of $a$. The \emph{type} of a tuple $(a^1,\ldots,a^n)$ of elements $a^1,\ldots,a^n\in \Xi_1\mult\cdots\mult\Xi_m$, denoted by $\typ(a^1,\ldots,a^n)$, is the $m$-tuple containing the types
    of $(a^1_i,\ldots,a^n_i)$ in $\Xi_i$ for each $1\leq i\leq m$. A function $f:\Xi_1\mult\cdots\mult\Xi_m\To \Omega$ is called \emph{canonical} iff it sends finite tuples of equal type in $\Xi_1\mult\cdots\mult\Xi_m$ to tuples of equal type in $\Omega$; that is, whenever $\typ(a^1,\ldots,a^n)=\typ(b^1,\ldots,b^n)$, then $\typ(f(a^1),\ldots,f(a^n))=\typ(f(b^1),\ldots,f(b^n))$. 
    For a relational structure $\Delta$ and elements $c_1,\ldots,c_n$ of $\Delta$, we write $(\Delta,c_1,\ldots,c_n)$ for the expansion of $\Delta$ by the constants $c_1,\ldots,c_n$. The structure $\Delta$ is \emph{ordered} iff it has a linear order among its relations.
     Now the following holds~\cite{BodPinTsa, BP-reductsRamsey}.
     
\begin{thm}[Bodirsky, Pinsker and Tsankov~\cite{BodPinTsa}]\label{thm:canonizing}
Let $\Delta$ be an ordered homogeneous Ramsey structure in a finite relational language, and let $\C\supseteq \Aut(\Delta)$ be a closed function clone. Then for all $f\in\C$ and all $c_1,\ldots,c_n\in\Delta$ there exists a function $g\in\C$ which is canonical as a function on $(\Delta,c_1,\ldots,c_n)$, and which agrees with $f$ on $\{c_1,\ldots,c_n\}$.
\end{thm}

Thus under these conditions on  $\Delta$, if there is a function $f$ in a polymorphism clone of a reduct of $\Delta$ which does something of interest (e.g., algorithmically) on a finite set $\{c_1,\ldots,c_n\}$, then there is also a canonical function in this clone which does the same.  Note that canonical functions on $(\Delta,c_1,\ldots,c_n)$
are finite objects in the following sense. Every canonical function $f:(\Delta,c_1,\ldots,c_n)^m\To (\Delta,c_1,\ldots,c_n)$ defines an $m$-ary function $T(f)$ on the types of $(\Delta,c_1,\ldots,c_n)$ in an obvious way, by the very definition of canonicity. Moreover, this type function $T(f)$ determines $f$ in the sense that if two canonical functions $f,g:(\Delta,c_1,\ldots,c_n)^m\To (\Delta,c_1,\ldots,c_n)$ have identical type functions $T(f)=T(g)$, then any closed function clone containing $\Aut(\Delta)$ contains $f$ iff it contains $g$. Since $(\Delta,c_1,\ldots,c_n)$ is homogeneous in a finite language, the type functions are finite objects: $T(f)$ is completely determined by its values on the types of tuples of length $q$, where $q$ is the maximal arity of a relation in $(\Delta,c_1,\ldots,c_n)$; moreover, there are only finitely many types of $q$-tuples since $(\Delta,c_1,\ldots,c_n)$ is $\omega$-categorical. As finite objects, these type functions can effectively be used in algorithms~\cite{BodPin-Schaefer-both}.

An example of an application of Theorem~\ref{thm:canonizing} is the following.

 \begin{thm}[Bodirsky, Pinsker and Tsankov~\cite{BodPinTsa}]\label{thm:canonical}
    Let $\Delta$ be an ordered homogeneous Ramsey structure in a finite relational language, and let $\Gamma$ be a reduct in a finite language. Then there exist $c_1,\ldots,c_n$, $m\geq 1$, and $m$-ary canonical functions $f_1,\ldots,f_k$ on $(\Delta,c_1,\ldots,c_n)$ such that for all reducts $\Gamma'$ we have that $\pol(\Gamma')\setminus\pol(\Gamma)$ is either empty or contains one of the functions $f_1,\ldots,f_k$.
\end{thm}

In words, under the above conditions the finite language reducts of $\Delta$ can be distinguished by functions which are canonical after adding finitely many constants to the language of $\Delta$. If we assume that that $\Delta$ is finitely bounded, which makes $\Delta$ in a way finitely representable, then Theorem~\ref{thm:canonical} can even be implemented in an algorithm. This yields the following effective variant of the theorem.

\begin{thm}[Bodirsky, Pinsker and Tsankov~\cite{BodPinTsa}]\label{thm:decidability}
Let $\Delta$ be an ordered homogeneous Ramsey structure which is finitely bounded, and let $\Gamma, \Gamma'$ be finite language reducts of $\Delta$. Then the problem whether or not $\Pol(\Gamma)\subseteq \Pol(\Gamma')$, where the relations of $\Gamma$ and $\Gamma'$ are given by quantifier-free formulas over $\Delta$, is decidable.
\end{thm}

This gives hope that tractability of CSPs of reducts of $\Delta$ is captured by the canonical functions in their polymorphism clones -- cf.~Section~\ref{sect:innovative}.

The modern proof of Theorem~\ref{thm:canonical} (yet unpublished but available on request) is based on a beautiful characterization of the Ramsey property for homogeneous structures which links Ramsey theory with topological dynamics~\cite{Topo-Dynamics}. A topological group $\G$ is called \emph{extremely amenable} iff whenever it acts continuously on a compact Hausdorff topological space $X$, then this action has a fixed point, i.e., there exists $x\in X$ such that $g(x)=x$ for all $g\in\G$. Let $\Delta$ be an ordered homogeneous relational structure. Then  $\Delta$ is Ramsey iff $\Aut(\Delta)$, viewed as an abstract topological group, is extremely amenable.

\section{Open Problems}\label{sect:innovative}

The research questions presented here are all related to Conjecture~\ref{conj:main} in one way or another; in fact, one can put them together so that they constitute a systematic program for proving the conjecture. After stating the questions, in Section~\ref{sect:conjecture}, I will discuss the questions in the context of the conjecture. I emphasize, however, that each of them has its own mathematical value independently of the truth of the conjecture.

The first set of questions concerns the connection between the algebraic and the topological structure of clones in the light of Theorem~\ref{thm:topobirk}. As for the link to constraint satisfaction, recall that for finite structures $\Gamma$, the complexity of $\Csp(\Gamma)$ only depends on the algebraic structure of $\pol(\Gamma)$, whereas in the $\omega$-categorical setting, one also has to take the topology on $\pol(\Gamma)$ into consideration (Theorem~\ref{thm:csp:topological}). A first question, which we shall then refine, is the following.

\begin{quest}\label{quest:topobirk}
 Are there conditions on oligomorphic algebras under which we can drop the continuity condition in Theorem~\ref{thm:topobirk}?
 \end{quest}

One such condition is to allow only very simple algebras $\B$. Of particular interest are, of course, continuous clone homomorphisms to {\bf 1} (i.e., term clones of algebras all of whose functions are projections), since they are the major source of hardness proofs. In fact, we believe that for certain structures they are the \emph{unique} source of hardness proofs: as already mentioned, the finite tractability conjecture states that if a finite structure $\Gamma$ satisfies a cosmetic condition, then $\Csp(\Gamma)$ is NP-hard if and only if there exists a homomorphism from $\pol(\Gamma)$ to {\bf 1}. That condition, which requires that $\pol(\Gamma)$ is \emph{idempotent}, i.e., all  $f\in\pol(\Gamma)$ satisfy the equation $f(x,\ldots,x)=x$, can always be assumed: every CSP of a finite structure is equivalent to a CSP of a finite structure with an idempotent polymorphism clone~\cite{JBK}. This fact can be derived in two steps: first, one shows that $\Gamma$ can be assumed to be a \emph{core}, i.e., all endomorphisms of $\Gamma$ are automorphisms. One then shows that it is possible to add finitely many constants to the language of $\Gamma$ without increasing the complexity of its CSP -- adding one constant for each element of the domain, this forces all polymorphisms to be idempotent. In the $\omega$-categorical setting, one cannot simply assume that $\pol(\Gamma)$ be idempotent; indeed, it would then certainly fail to be oligomorphic, containing no unary functions at all except the identity. However, it is possible to perform an analog  of the first step, and assume that $\Gamma$ is a \emph{model-complete core}, meaning that $\Aut(\Gamma)$ is (topologically) dense in the endomorphism monoid of $\Gamma$ (i.e., every endomorphism \emph{locally} looks like an automorphism)~\cite{Bodirsky-HDR}. Moreover, as in the finite case, adding finitely many constants to the language of $\Gamma$ then does not increase the complexity of its CSP~\cite{Bodirsky-HDR}. Until now, similarly to the situation for finite templates, we do not know of a reduct of a finitely bounded homogeneous structure which is a model-complete core with an NP-hard CSP, but which has no continuous homomorphism to {\bf 1} after adding finitely many constants to the language.

Back to the continuity condition, we therefore ask the following.

\begin{quest}\label{quest:homotoone}
If a closed oligomorphic clone has a clone homomorphism to {\bf 1} (i.e., it satisfies no non-trivial equations), does it always have a continuous clone homomorphism to {\bf 1}? If not, are there further conditions on the clone (model complete core etc.) which imply a positive answer?
\end{quest}

There exists considerable literature 
about automorphism groups of $\omega$-categorical structures $\Gamma$ which are \emph{reconstrucible}, i.e., where the topology
on $\Aut(\Gamma)$
is uniquely determined by the
algebraic group structure; this is for instance the case when  
$\Aut(\Gamma)$ has the so-called \emph{small index property},
that is, all subgroups of countable index 
are open. The small index property has for instance been shown for $\Aut(\mathbb N,=)$~\cite{DixonNeumannThomas}; for $\Aut({\mathbb Q};<)$ and for the automorphism group of the atomless Boolean algebra~\cite{Truss}; the automorphism group of the random graph~\cite{HodgesHodkinsonLascarShelah}; for all $\omega$-categorical $\omega$-stable structures~\cite{HodgesHodkinsonLascarShelah}; for the automorphism groups of the Henson graphs~\cite{Herwig98}. The notion of reconstruction makes perfect sense for function clones,  and is of importance for our purposes. Call a closed function clone \emph{reconstructible} iff all isomorphisms with other closed function clones are homeomorphisms. Recent research has shown that for some homogeneous $\omega$-categorical structures with a reconstructible automorphism group, the reconstructability carries over to the polymorphism clone of the structure~\cite{BPP-reconstructing}.

\begin{quest}\label{quest:automaticContinuity}
Let $\C$ be an oligomorphic polymorphism clone whose group of invertible unary functions is reconstructible. When can we conclude that $\C$ is reconstructible  as well?
\end{quest}

We remark that there exists an example of two $\omega$-categorical structures whose automorphism groups are isomorphic as groups but not as topological groups~\cite{EvansHewitt}, and that this example has recently been expanded to polymorphism clones by David Evans in a yet unpublished note.

 It is well-known that every Baire measurable  homomorphism between Polish groups is continuous (see e.g.~\cite{Kechris}).
So let us remark that there exists a model of ZF+DC where every set is Baire measurable~\cite{Shelah84}. 
For the structures $\Gamma$ that we need to model computational problems as $\Csp(\Gamma)$
it therefore seems fair to assume that the abstract algebraic structure of $\Aut(\Gamma)$ always determines its topological structure; consistency of this statement with ZF has already been observed in~\cite{Lascar}. Hence, one could hope to find a model of ZF in which polymorphism clones of $\omega$-categorical structures are reconstructible, or in which all homomorphisms of such clones to {\bf 1} are continuous.

\begin{quest}\label{quest:ZF}
Do oligomorphic polymorphism clones have reconstruction in an appropriate model of ZF? Are all homomorphisms from oligomorphic polymorphism clones to {\bf 1}  continuous in an appropriate model of ZF?
\end{quest}

The concept of a \emph{topological clone} appeared as a necessity for formulating Theorem~\ref{thm:topobirk}; it is indeed natural considering the importance of abstract clones (known in disguise as \emph{varieties}) for universal algebra and the natural presence of topological groups in various fields of mathematics, in particular topological dynamics. It is known that the closed permutation groups on a countable set  
are precisely those topological groups that are Polish
and have a left-invariant ultrametric~\cite{BeckerKechris}.

\begin{quest}\label{quest:topoclones}
Which topological clones appear as closed function clones on a countably infinite set?
\end{quest}

We now turn to the study of polymorphism clones of reducts of homogeneous Ramsey structures. Here, the approach via canonical functions, based on Theorem~\ref{thm:canonizing} and the idea that we keep sufficient information about a clone when we add a sufficiently large finite number of constants and then only consider its canonical functions, has proven extremely fruitful. For example, this was the strategy in the Graph-SAT dichotomy classification~\cite{BodPin-Schaefer-both}, and in many other applications~\cite{Poset-Reducts,42,BodPinTsa,BP-reductsRamsey,RandomMinOps}; confer also  Theorems~\ref{thm:canonical} and~\ref{thm:decidability}. Generalizing the Graph-SAT strategy, we arrive at the following ideas.

In the following, let $\Gamma$ be a reduct of a finitely bounded ordered homogeneous Ramsey structure $\Delta$, and let $c_1,\ldots,c_n\in\Delta$. Then the set of finitary canonical functions on $(\Delta,c_1,\ldots,c_n)$ forms a closed function clone, and hence so does the intersection of this clone with $\pol(\Gamma)$, which we call the \emph{canonical fragment} of $\pol(\Gamma)$ with respect to $c_1,\ldots,c_n$. By Theorem~\ref{thm:canonical}, this canonical fragment still contains considerable information about $\pol(\Gamma)$. Recall that its functions define functions on the types of $(\Delta,c_1,\ldots,c_n)$; in fact, these ``type functions'' form a clone on a finite set. We call this clone the \emph{type clone} of $\pol(\Gamma)$ with respect to $c_1,\ldots,c_n$, and denote it by $T_{c_1,\ldots,c_n}(\pol(\Gamma))$. There are infinitely many choices for $c_1,\ldots,c_n$, but up to type equivalence, only finitely many for each $n$ since $\Gamma$ is $\omega$-categorical~\cite{HodgesLong}. In the Graph-SAT dichotomy, these type clones happened to capture the computational complexity of $\Csp(\Gamma)$~\cite{BodPin-Schaefer-both}.

\begin{quest}\label{quest:types:general}
Does the complexity of $\Csp(\Gamma)$ only depend on the algebraic structure of its type clones?
\end{quest}

More precisely, in the Graph-SAT classification it turned out that when $\Csp(\Gamma)$ is tractable, then this fact was captured by some canonical polymorphism, which provided the algorithm; in particular, the answer to the following question was positive. It is nourished by the belief that if $\pol(\Gamma)$ contains a function which implies tractability, and therefore is of use in some algorithm, then this function can be ``canonized'' and therefore appears in some type clone (cf.~Theorem~\ref{thm:canonizing}).

\begin{quest}\label{quest:types:tractability:down}
If $\Gamma$ is tractable, are there necessarily $c_1,\ldots,c_n\in\Delta$ such that $T_{c_1,\ldots,c_n}(\pol(\Gamma))$ corresponds to a tractable CSP?
\end{quest}

The following question asks about the converse. The intuition behind it, again true in the Graph-SAT case, is that algorithms for the type clone can be ``lifted'' back to the original clone.

\begin{quest}\label{quest:types:tractability:up}
If there exist $c_1,\ldots,c_n\in\Delta$ such that $T_{c_1,\ldots,c_n}(\pol(\Gamma))$ is tractable, is $\pol(\Gamma)$ tractable?
\end{quest}

The following can be seen as a complexity-free variant of the preceding two questions. It is known that if $T_{c_1,\ldots,c_n}(\pol(\Gamma))$ satisfies non-trivial equations for some $c_1,\ldots,c_n\in\Delta$, then so does $\pol(\Gamma)$~\cite{BPP-projective}. Therefore, if some type clone corresponds to a tractable, and if the finite tractability conjecture holds, then $\pol(\Gamma)$ does not have a homomorphism to {\bf 1} (and hence, if NP-hard, would have to have another source of hardness). We do not know if the converse is true as well.

\begin{quest}\label{quest:equations}
If $\pol(\Gamma)$ does not have a homomorphism to $\mathbf 1$, are there $c_1,\ldots,c_n\in\Delta$ such that $T_{c_1,\ldots,c_n}(\pol(\Gamma))$ has no homomorphism to $\mathbf 1$?
\end{quest}

We often think of the Ramsey property as an additional property of Fra\"{i}ss\'{e} classes (and indeed, most homogeneous structures are not Ramsey); we do require the property for our methods. However, we do not require the reducts $\Gamma$, but only some ``base structure'' $\Delta$ in which they are definable, to be Ramsey. Hence, if the answer to the following question were positive, then our methods would work for all homogeneous structures.

\begin{quest}\label{quest:RamseyExpansion}
Can every finitely bounded Fra\"{i}ss\'{e} class be extended by finitely many relations to a finitely bounded Fra\"{i}ss\'{e} class 
 which is in addition Ramsey?
\end{quest}
We remark that this question has recently received considerable attention, and in particular the answer is positive for all Fra\"{i}ss\'{e} classes of digraphs -- see~\cite{LafJasTheWoo} for further references.

Let us turn to concrete classes of CSPs. Studying those has its own interest, just like the dichotomies of Graph-SAT problems~\cite{BodPin-Schaefer-both} and temporal constraints~\cite{tcsps-journal, tcsps}; moreover they provide sources of examples for the general questions.

The class of finite partial orders is one of the most natural Fra\"{i}ss\'{e} classes, and the answer to the following question would subsume some older results in theoretical computer science, e.g. in~\cite{BroxvallJonsson}, and the classification~\cite{tcsps-journal}. The basis for a successful complexity classification of Poset-SAT problems has been established very recently in the form of the classification of the closed supergroups of the automorphism group of the random partial order~\cite{Poset-Reducts}.

\begin{quest}\label{quest:complexity:poset}
Classify the complexity of $\Csp(\Gamma)$, for all finite language reducts $\Gamma$ of the random partial order. In other words, classify the complexity of Poset-SAT problems.
\end{quest}

The solution to the following would subsume a considerable amount of results in the literature, e.g., completely the papers~\cite{RCC5JD}, \cite{BodHils}, and some results in~\cite{BroxvallJonsson},~\cite{MarriottOdersky}. It would moreover require the extension of our methods to \emph{functional} signatures, a venture interesting in itself.

\begin{quest}\label{quest:complexity:BA}
Classify the complexity of $\Csp(\Gamma)$, for all finite language reducts $\Gamma$ of the atomless (= random) Boolean algebra.
\end{quest}

\section{The Infinite Tractability Conjecture}\label{sect:conjecture}

Each of the above questions has its own interest for the understanding of oligomorphic function clones, oligomorphic algebras, and their connections with constraint satisfaction. However, these questions really are part of a bigger program around Conjecture~\ref{conj:main}. To make this connection with Conjecture~\ref{conj:main} more evident, let me show an example of how a proof of the conjecture could look like. Let $\Gamma$ be a finite language reduct of a finitely bounded homogeneous structure $\Delta$.

\begin{itemize}
\item Assume that $\Pol(\Gamma)$ does not have a continuous clone homomorphism to {\bf 1} (otherwise, $\Csp(\Gamma)$ is NP-hard by Theorem~\ref{thm:csp:topological}).
\item If the answer to Question~\ref{quest:homotoone} is positive (possibly in an appropriate model of ZF, cf.~Question~\ref{quest:ZF}), then $\pol(\Gamma)$ satisfies non-trivial equations.
\item Assuming that Question~\ref{quest:RamseyExpansion} has a positive answer, we can assume that $\Delta$ is Ramsey.
\item If Question~\ref{quest:equations} has a positive answer, then some type clone $T_{c_1,\ldots,c_n}(\pol(\Gamma))$ satisfies non-trivial equations.
\item Assuming the tractability conjecture for finite templates, $T_{c_1,\ldots,c_n}(\pol(\Gamma))$ corresponds to a tractable CSP.
\item Assuming a positive answer to Question~\ref{quest:types:tractability:up}, $\Csp(\Gamma)$ is then tractable.
\end{itemize}

\bibliographystyle{plain}
\bibliography{doccourse}

\end{document}

%% file: theoremstyle.tex
\theoremstyle{plain}
    \newtheorem{thm}{Theorem}

    \newtheorem{quest}[thm]{Question}
    
    \newtheorem{conj}[thm]{Conjecture}

\theoremstyle{definition}

    \newtheorem{nota}[thm]{Notation}
\theoremstyle{remark}

%% file: commands.tex
\newcommand{\ignore}[1]{}

\newcommand{\HSPfin}{\HSP^{\fin}}

\newcommand{\cC}{\ensuremath{\mathcal{C}}}
\newcommand{\cV}{\ensuremath{\mathcal{V}}}

\DeclareMathOperator{\fin}{fin}

\DeclareMathOperator{\pol}{Pol}

\DeclareMathOperator{\Clo}{Clo}
\DeclareMathOperator{\HSP}{HSP}
\DeclareMathOperator{\HHH}{H}
\DeclareMathOperator{\PPP}{P}
\DeclareMathOperator{\PPPfin}{P^{\fin}}
\DeclareMathOperator{\SSS}{S}

\DeclareMathOperator{\typ}{tp}

 \DeclareMathOperator{\Csp}{CSP}

\DeclareMathOperator{\Aut}{Aut}
\DeclareMathOperator{\Pol}{Pol}

\newcommand{\C}{{\mathscr C}}

\newcommand{\F}{{\mathscr F}}

\newcommand{\A}{{\mathfrak A}}
\newcommand{\B}{{\mathfrak B}}
\newcommand{\D}{{\mathscr D}}

\newcommand{\G}{{\mathscr G}}

\newcommand{\To}{\rightarrow}
\newcommand{\mult}{\times}

\newcommand{\cproblem}[3]{
\vspace{.2cm}
\noindent {\bf #1} \\
INSTANCE: #2 \\
QUESTION: #3 \\}

\def\ocirc#1{\ifmmode\setbox0=\hbox{$#1$}\dimen0=\ht0
    \advance\dimen0 by1pt\rlap{\hbox to\wd0{\hss\raise\dimen0
    \hbox{\hskip.2em$\scriptscriptstyle\circ$}\hss}}#1\else
    {\accent"17 #1}\fi}